\documentclass[12pt,reqno]{article}
\usepackage{amssymb}
\usepackage{amsthm}
\usepackage{amsmath}
\newtheorem{theorem}{Theorem}[section]

\newtheorem{example}[theorem]{Example}

\newtheorem{remark}[theorem]{Remark}

\def\be{\begin{equation}}
\def\ee{\end{equation}}
\def\bea{\begin{eqnarray}}
\def\eea{\end{eqnarray}}
\begin{document}
\begin{center} \Large{\bf Invariant characterization of scalar third-order ODEs that admit the maximal point symmetry Lie algebra}
\end{center}
\medskip
\begin{center}
Ahmad Y. Al-Dweik$^*$, M. T. Mustafa$^{**}$ and
F. M. Mahomed$^{***,}$\\

{$^*$Department of Mathematics \& Statistics, King Fahd University
of Petroleum and Minerals, Dhahran 31261, Saudi Arabia}\\
{$^{**}$Department of Mathematics, Statistics and Physics, Qatar
University, Doha, 2713, State of Qatar}\\
{$^{***}$School of Computer Science and Applied Mathematics,
DST-NRF Centre of Excellence in Mathematical and Statistical
Sciences,
University of the Witwatersrand, Johannesburg, Wits 2050,  South Africa\\
}

aydweik@kfupm.edu.sa, tahir.mustafa@qu.edu.qa and
Fazal.Mahomed@wits.ac.za
\end{center}
\begin{abstract}
The Cartan equivalence method is utilized to deduce an invariant
characterization of the scalar third-order ordinary differential
equation $u'''=f(x,u,u',u'')$ which admits the maximal
seven-dimensional point symmetry Lie algebra. The method provides
auxiliary functions which can be used to efficiently  obtain the
point transformation that does the reduction to the simplest
linear equation $\bar{u}'''=0$. Moreover, examples are given to
illustrate the method.
\end{abstract}
\newpage
\section{Introduction}
Lie \cite{lie,lie1} was the first to determine practical and
algebraic  criteria for linearization for scalar second-order
ordinary differential equations (ODEs) via invertible maps. The
reduction to the simplest form was for equations which
possessed the maximal symmetry algebra of dimension eight.
Tress\'e then showed that any scalar second-order ODE
 is linearizable by means of
a point transformation if and only if the two Tress\'e \cite{tre}
relative invariants both vanish for the ODE.

Cartan \cite{car} solved, inter alia, the linearzation problem for
ODEs using his approach now called the Cartan equivalent method
(see the insightful contributions \cite{Gardner1989,Olver1995} on
this). In fact it was Grissom et. al. \cite{gri} who first used
the Cartan equivalence method on the Lie problem to obtain the Lie
invariant criteria for linearization for scalar second-order ODEs
in a geometric manner.

The Lie algebraic criteria for linearization by means of point
transformations for scalar $n$th-order ($n>2$) ODEs were uncovered
in \cite{mah4}. Three canonical forms arose for scalar linear
third-order ODEs. The maximal Lie algebra case for such ODEs is
dimension seven and corresponds to $u'''=0$.

The Laguerre-Forsyth (see \cite{mah4}) canonical form for scalar
linear third-order ODEs is
\begin{equation}\label{e3}
u'''+a(x)u=0.
\end{equation} If $a\ne0$, then (\ref{e3}) has a five- or four-dimensional symmetry
algebra.

Chern \cite{che} applied the Cartan equivalence method to solve
the linearization problem of scalar third-order ODEs via contact
transformations. He obtained conditions of equivalence for
(\ref{e3}) in the cases $a=0$ and $a=1$. Neut and Petitot
\cite{neu} also studied conditions on equivalence to (\ref{e3}).
Grebot \cite{gre} looked at linearization of third-order ODEs via
fibre preserving point transformations. In recent work, Ibragimov
and Meleshko \cite{ibr} studied the linearization problem using a
direct Lie approach for such third-order ODEs via both point and
contact transformations. They also gave conditions on the
linearizing maps.

We investigate herein the linearization problem via invertible
point transformation for scalar third-order ODEs which admit the
maximal dimension point symmetry algebra via the Cartan approach.
We also provide compact criteria in terms of $f$. We previously in
\cite{DweikTahirFazal} considered the submaximal five-dimensional
case for linearization invoking this method to provide compact
criteria in terms of $f$. The case of seven point symmetries was
also studied by Al-Dweik \cite{Dweik2016}. Here we refine this and
also determine conditions on how to find the invertible maps using
auxiliary functions.

We mention here that invariant characterization of
$u'''=f(x,u,u',u'')$ that admits the maximal point symmetry
algebra was obtained in terms of the function $f$ in the following
theorem. We denote $u', u''$ by $p, q$, respectively in the
following and thereafter.
\begin{theorem}\cite{Dweik2016}
The necessary and sufficient conditions for equivalence of a
third-order equation $ u''' = f(x,u,u',u'')$ to the canonical form
$u''' = 0$ with seven symmetries via {\rm point transformation}
are the identical vanishing of the system of relative invariants
\begin{equation}\label{c1}
\begin{array}{ll}
{f_{q,q,q} }  \\
{f_{q,q} ^2  + 6\,f_{p,q,q} }  \\
4\,f_q ^3  + 18\,f_q \left( {f_p  - \dot{D}_x f_q} \right) + 9\,\dot{D}_x^2 f_q  + 54\,f_u  - 27\,\dot{D}_x f_p\\
f_{q,q} \left( {f_q ^2+9\,f_p  - 3\,\dot{D}_x f_q}\right) - 9\,f_{p,p}  + 18\,f_{u,q}  - 6\,f_q f_{p,q},\\
\end{array}
\end{equation}
where $ \dot{D}_x = \frac{\partial } {{\partial x}} +
p\frac{\partial }{{\partial u}}+q\frac{\partial} {{\partial
p}}+f\frac{\partial} {{\partial q}}$.
\end{theorem}
\section{Application of Cartan's equivalence method for third-order ODEs with seven point symmetries}
Let $x:=(x,u,p=u',q=u'')\in \mathbb{R}^4$ as usual be local coordinates
of $J^2$, the space of the second order jets. In local
coordinates, the equivalence of two third-order ODEs
\begin{equation}\label{b0}
\begin{array}{cc}
u'''=f(x,u,u',u''), & \bar{u}'''=\bar{f}(\bar{x},\bar{u},\bar{u}',\bar{u}''),\\
\end{array}
\end{equation}
under an invertible transformation
\begin{equation}\label{ccc}
\bar{x}=\phi \left( x,u \right),~\bar{u} =\psi \left( x,u  \right),~~~\phi_x\psi_u   -  \phi_u \psi_x\neq0,\\
\end{equation}
can be expressed as the local equivalence problem for the $G$-structure
\begin{equation}\label{b2}
\Phi^*\left(%
\begin{array}{c}
  \bar{\omega}^1 \\
  \bar{\omega}^2 \\
  \bar{\omega}^3 \\
  \bar{\omega}^4 \\
\end{array}%
\right)=\left(%
\begin{array}{cccc}
  a_1 & 0 & 0 & 0 \\
  a_2 & a_3 & 0 & 0 \\
  a_4 & a_5 & a_6 & 0 \\
  a_7 & 0 & 0 & a_8 \\
\end{array}%
\right) \left(%
\begin{array}{c}
  \omega^1 \\
  \omega^2 \\
  \omega^3 \\
  \omega^4 \\
\end{array}%
\right),
\end{equation}
where
\begin{equation}\label{b1}
\begin{array}{llll}
\bar{\omega}^1=d\bar{u}-\bar{p} d\bar{x}, & \bar{\omega}^2=d\bar{p}-\bar{q} d \bar{x}, & \bar{\omega}^3=d\bar{q}-\bar{f} d \bar{x}, & \bar{\omega}^4= d \bar{x},\\
\omega^1=du-p d x, & \omega^2=dp-q d x, & \omega^3=dq-f d x, & \omega^4= d x.\\
\end{array}
\end{equation}
One can evaluate the functions $a_i=a_i(x,u,p,q), i=1..8,$
explicitly as $a_1=\frac{\phi_x\psi_u   - \phi_u
\psi_x}{D_x \phi},\\ a_2(x,u,p,q)=\frac{D_x a_1}{D_x \phi},
a_3(x,u,p)=\frac{a_1}{D_x \phi}$.

Now, one can define $\theta$ to be the lifted coframe with an
eight-dimensional group $G$
\begin{equation}\label{b21}
\left(%
\begin{array}{c}
  \theta^1 \\
  \theta^2 \\
  \theta^3 \\
  \theta^4 \\
\end{array}%
\right)=\left(%
\begin{array}{cccc}
  a_1 & 0 & 0 & 0 \\
  a_2 & a_3 & 0 & 0 \\
  a_4 & a_5 & a_6 & 0 \\
  a_7 & 0 & 0 & a_8 \\
\end{array}%
\right) \left(%
\begin{array}{c}
  \omega^1 \\
  \omega^2 \\
  \omega^3 \\
  \omega^4 \\
\end{array}%
\right).
\end{equation}
The application of Cartan's method to this equivalence problem leads
to an ${e}$-structure, which is invariantly associated to the
given equation. We mimic the sequence of calculations given in our work \cite{DweikTahirFazal}.
This is needed in what follows.

We have that the first structure equation is
\begin{equation}\label{b3}
d\left(%
\begin{array}{c}
  \theta^1 \\
  \theta^2 \\
  \theta^3 \\
  \theta^4 \\
\end{array}%
\right)=\left(%
\begin{array}{cccc}
  \alpha_1 & 0 & 0 & 0 \\
  \alpha_2 & \alpha_3 & 0 & 0 \\
  \alpha_4 & \alpha_5 & \alpha_6 & 0 \\
  \alpha_7 & 0 & 0 & \alpha_8 \\
\end{array}%
\right)\wedge \left(%
\begin{array}{c}
  \theta^1 \\
  \theta^2 \\
  \theta^3 \\
  \theta^4 \\
\end{array}%
\right)+
\left(%
\begin{array}{c}
  T^1_{24}~\theta^2 \wedge \theta^4 \\
  T^2_{34}~\theta^3 \wedge \theta^4  \\
  0 \\
  0 \\
\end{array}%
\right)\\
\end{equation}
The infinitesimal action on the torsion is of the form
\begin{equation}\label{b4}
\left.
\begin{array}{cccc}
d~T^1_{24}\equiv(\alpha_1-\alpha_3-\alpha_8)T^1_{24}\\
d~T^2_{34}\equiv(\alpha_3-\alpha_6-\alpha_8)T^2_{34}\\
\end{array}
\right\}
~\textrm{mod}~(\theta^1,\theta^2,\theta^3,\theta^4)\\
\end{equation}
and a parametric evaluation gives $T^1_{24}=-\frac{a_1}{a_3
a_8}\ne0$ and $T^2_{34}=-\frac{a_3}{a_6 a_8}\ne0$. One normalizes
the torsion by
\begin{equation}\label{b5}
\begin{array}{cc}
T^1_{24}=-1,&T^2_{34}=-1.\\
\end{array}
\end{equation}
This results in the principal components as follows
\begin{equation}\label{b6}
\begin{array}{cc}
\alpha_6=2\alpha_3-\alpha_1, &\alpha_8=\alpha_1-\alpha_3.\\
\end{array}
\end{equation}
The normalizations enforce relations on the group $G$ as
\begin{equation}\label{b61}
\begin{array}{cc}
a_6=\frac{a_3^2}{a_1},&a_8=\frac{a_1}{a_3}.\\
\end{array}
\end{equation}
The {\it first-order} normalizations gives rise to an adapted coframe
with the {\it six-dimensional group} $G_1$
\begin{equation}\label{b62}
\left(%
\begin{array}{c}
  \theta^1 \\
  \theta^2 \\
  \theta^3 \\
  \theta^4 \\
\end{array}%
\right)=\left(%
\begin{array}{cccc}
  a_1 & 0 & 0 & 0 \\
  a_2 & a_3 & 0 & 0 \\
  a_4 & a_5 & \frac{a_3^2}{a_1} & 0 \\
  a_7 & 0 & 0 & \frac{a_1}{a_3} \\
\end{array}%
\right) \left(%
\begin{array}{c}
  \omega^1 \\
  \omega^2 \\
  \omega^3 \\
  \omega^4 \\
\end{array}%
\right).
\end{equation}
This yields the structure equation
\begin{equation}\label{b7}
d\left(%
\begin{array}{c}
  \theta^1 \\
  \theta^2 \\
  \theta^3 \\
  \theta^4 \\
\end{array}%
\right)=\left(%
\begin{array}{cccc}
  \alpha_1 & 0 & 0 & 0 \\
  \alpha_2 & \alpha_3 & 0 & 0 \\
  \alpha_4 & \alpha_5 & 2\alpha_3-\alpha_1 & 0 \\
  \alpha_7 & 0 & 0 & \alpha_1-\alpha_3 \\
\end{array}%
\right)\wedge \left(%
\begin{array}{c}
  \theta^1 \\
  \theta^2 \\
  \theta^3 \\
  \theta^4 \\
\end{array}%
\right)+
\left(%
\begin{array}{c}
  -\theta^2 \wedge \theta^4 \\
 -\theta^3 \wedge \theta^4  \\
  T^3_{34}~\theta^3 \wedge \theta^4  \\
  0  \\
\end{array}%
\right)\\
\end{equation}
The infinitesimal action on the torsion is
\begin{equation}\label{b8}
\begin{array}{cccc}
d~T^3_{34}\equiv (\alpha_3-\alpha_1)T^3_{34}+3(\alpha_2-\alpha_5)\\
\end{array}
~\textrm{mod}~(\theta^1,\theta^2,\theta^3,\theta^4)\\
\end{equation}
and one can translate $T^3_{34}$ to vanish as
\begin{equation}\label{b9}
\begin{array}{cc}
T^3_{34}=0.\\
\end{array}
\end{equation}
This gives the principal components as
\begin{equation}\label{b91}
\begin{array}{cc}
\alpha_5=\alpha_2.\\
\end{array}
\end{equation}
The normalizations force relations on the group $G_1$ to be
\begin{equation}\label{b92}
\begin{array}{c}
a_5=\frac{a_2a_3}{a_1}-\frac{a_3^2}{3a_1} s_1,\\
\end{array}
\end{equation}
where $s_1=f_q.$

The {\it second-order} normalizations results in  an adapted coframe
with the {\it five-dimensional group} $G_2$
\begin{equation}\label{b93}
\left(%
\begin{array}{c}
  \theta^1 \\
  \theta^2 \\
  \theta^3 \\
  \theta^4 \\
\end{array}%
\right)=\left(%
\begin{array}{cccc}
  a_1 & 0 & 0 & 0 \\
  a_2 & a_3 & 0 & 0 \\
  a_4 & \frac{a_2a_3}{a_1}-\frac{a_3^2}{3a_1} s_1 & \frac{a_3^2}{a_1} & 0 \\
  a_7 & 0 & 0 & \frac{a_1}{a_3} \\
\end{array}%
\right) \left(%
\begin{array}{c}
  \omega^1 \\
  \omega^2 \\
  \omega^3 \\
  \omega^4 \\
\end{array}%
\right).
\end{equation}
This yields the structure equation
\begin{equation}\label{b10}
d\left(%
\begin{array}{c}
  \theta^1 \\
  \theta^2 \\
  \theta^3 \\
  \theta^4 \\
\end{array}%
\right)=\left(%
\begin{array}{cccc}
  \alpha_1 & 0 & 0 & 0 \\
  \alpha_2 & \alpha_3 & 0 & 0 \\
  \alpha_4 & \alpha_2 & 2\alpha_3-\alpha_1 & 0 \\
  \alpha_7 & 0 & 0 & \alpha_1-\alpha_3 \\
\end{array}%
\right)\wedge \left(%
\begin{array}{c}
  \theta^1 \\
  \theta^2 \\
  \theta^3 \\
  \theta^4 \\
\end{array}%
\right)+
\left(%
\begin{array}{c}
  -\theta^2 \wedge \theta^4 \\
 -\theta^3 \wedge \theta^4  \\
  T^3_{24}~\theta^2 \wedge \theta^4  \\
  T^4_{24}~\theta^2 \wedge \theta^4  \\
\end{array}%
\right)\\
\end{equation}
The infinitesimal action on the torsion now is
\begin{equation}\label{b11}
\left.
\begin{array}{ll}
d~T^3_{24}\equiv 2(\alpha_3-\alpha_1)T^3_{24}-2\alpha_4\\
d~T^4_{24}\equiv-\alpha_3~T^4_{24}-\alpha_7\\
\end{array}
\right\}
~\textrm{mod}~(\theta^1,\theta^2,\theta^3,\theta^4)\\
\end{equation}
and we can translate $T^3_{24}$ and $T^4_{24}$  to zero
\begin{equation}\label{b12}
\begin{array}{cc}
T^3_{24}=0,&T^4_{24}=0.\\
\end{array}
\end{equation}
This gives rise to the principal components
\begin{equation}\label{b121}
\begin{array}{cc}
\alpha_4=0,&\alpha_7=0.\\
\end{array}
\end{equation}
The normalizations force relations on the group $G_2$ as
\begin{equation}\label{b122}
\begin{array}{cc}
a_4=\frac{a^2_2}{2a_1}-\frac{a^2_3}{18 a_1}s_2,& a_7=\frac{a_1}{6a_3}s_3,\\
\end{array}
\end{equation}
where $s_2= 2f_q ^2 +9\,f_p   - 3\,D_x f_q ,~s_3=f_{q,q}$.

The {\it third-order} normalizations gives rise to  an adapted coframe
with the {\it three-dimensional group} $G_3$
\begin{equation}\label{b123}
\left(%
\begin{array}{c}
  \theta^1 \\
  \theta^2 \\
  \theta^3 \\
  \theta^4 \\
\end{array}%
\right)=\left(%
\begin{array}{cccc}
  a_1 & 0 & 0 & 0 \\
  a_2 & a_3 & 0 & 0 \\
  \frac{a^2_2}{2a_1}-\frac{a^2_3}{18 a_1}s_2 & \frac{a_2a_3}{a_1}-\frac{a_3^2}{3a_1} s_1 & \frac{a_3^2}{a_1} & 0 \\
  \frac{a_1}{6a_3}s_3 & 0 & 0 & \frac{a_1}{a_3} \\
\end{array}%
\right) \left(%
\begin{array}{c}
  \omega^1 \\
  \omega^2 \\
  \omega^3 \\
  \omega^4 \\
\end{array}%
\right).
\end{equation}
This results in the structure equation
\begin{equation}\label{b13}
d\left(%
\begin{array}{c}
  \theta^1 \\
  \theta^2 \\
  \theta^3 \\
  \theta^4 \\
\end{array}%
\right)=\left(%
\begin{array}{cccc}
  \alpha_1 & 0 & 0 & 0 \\
  \alpha_2 & \alpha_3 & 0 & 0 \\
  0 & \alpha_2 & 2\alpha_3-\alpha_1 & 0 \\
  0 & 0 & 0 & \alpha_1-\alpha_3 \\
\end{array}%
\right)\wedge \left(%
\begin{array}{c}
  \theta^1 \\
  \theta^2 \\
  \theta^3 \\
  \theta^4 \\
\end{array}%
\right)+
\left(%
\begin{array}{c}
  -\theta^2 \wedge \theta^4 \\
 -\theta^3 \wedge \theta^4  \\
  T^3_{14}~\theta^1 \wedge \theta^4  \\
  T^4_{12}~\theta^1 \wedge \theta^2+T^4_{13}~\theta^1 \wedge \theta^3  \\
\end{array}%
\right)\\
\end{equation}
The infinitesimal action on the torsion is of the form
\begin{equation}\label{b14}
\left.
\begin{array}{cc}
d~T^3_{14}\equiv &-3(\alpha_1-\alpha_3) T^3_{14}\\
d~T^4_{12}\equiv &-2 \alpha_3 T^4_{12}-\alpha_2 T^4_{13} \\
d~T^4_{13}\equiv &(\alpha_1-3\alpha_3) T^4_{13}\\
\end{array}
\right\}
~\textrm{mod}~(\theta^1,\theta^2,\theta^3,\theta^4)\\
\end{equation}
and one has a bifurcation in the flowchart which depends upon
whether $T^3_{14}$, $T^4_{12}$ and $T^4_{13}$ are zero. A
parametric calculation results in
\begin{equation}\label{b15}
\begin{array}{ll}
T^3_{14}=&-\frac{1}{54}\frac{a_3^3 I_3}{a_1^3},\\
T^4_{13}=& -\frac{1}{6}\frac{a_1 I_1}{a_3^3}, \hfill \\
T^4_{12}=& -\frac{1}{36}\frac{I_2}{a_3^2}~~\textrm{mod}~~T^4_{13}, \hfill \\
\end{array}
\end{equation}
where
\begin{equation}\label{b16}
\begin{array}{ll}
  I_1  &={s_3}_q= f_{q,q,q},\\
  I_2  &={s^2_3}+6{s_3}_p =f_{q,q} ^2  + 6\,f_{p,q,q},  \hfill \\
  I_3  &=2s_1s_2-3D_x s_2+54f_u\\
       &= {4\,f_q ^3  + 18\,f_q \left( {f_p  - D_x f_q } \right) + 9\,D_x^2 f_q  - 27\,D_x f_p + 54\,f_u }.  \hfill \\
\end{array}
\end{equation}
It is well-known that third-order ODEs with seven point
symmetries have the canonical form $u''' =0$. It should be noted
here that the relative invariants $I_1=I_2=I_3=0$  for the
canonical form $u''' = 0$. Therefore, we choose the following
branch.
\section*{Branch 1. $I_1=I_2=I_3=0$.}
In this branch, there is no more unabsorbable torsion left, so the
remaining group variables $a_1, a_2$ and $a_3$ cannot be
normalized. Moreover, $\alpha_1, \alpha_2$ and $\alpha_3$ are
uniquely defined, so the problem is determinant. This results in
the following $e$-structure on the seven-dimensional prolonged space
$M^{(1)}=M \times G_3$ which consists of the original lifted coframe
\begin{equation}\label{b36}
\begin{array}{l}
  \theta^1=a_1 \omega^1,\\
  \theta^2=a_2 \omega^1+a_3 \omega^2,\\
  \theta^3=\left(\frac{a^2_2}{2a_1}-\frac{a^2_3}{18 a_1}s_2\right)\omega^1+\left(\frac{a_2a_3}{a_1}-\frac{a_3^2}{3a_1} s_1\right)\omega^2+\frac{a_3^2}{a_1}\omega^3,\\
  \theta^4= \frac{a_1}{6a_3}s_3~\omega^1+\frac{a_1}{a_3}~\omega^4,\\
 \end{array}
\end{equation}
together with the modified Maurer-Cartan forms
\begin{equation}
  \left.
  \begin{array}{l}
  \alpha_1=\frac{da_1}{a_1}-\frac{a_2}{a_3}~\omega^4,\\
  \alpha_2=\frac{da_2}{a_1}-\frac{a_2da_3}{a_1a_3}+\frac{1}{18}\left(\frac{a^2_3s_2-9a^2_2-6a_2a_3s_1}{a_1a_3}\right)~\omega^4,\\
  \alpha_3=\frac{da_3}{a_3}+\frac{1}{3}s_1~\omega^4.\\
 \end{array}
  \right\}
~\textrm{mod}~(\omega^1,\omega^2,\omega^3)\\
\end{equation}
This results in the structure equations
\begin{equation}\label{b35}
\begin{array}{l}
  d\theta^1=-\theta^1 \wedge \alpha_1-\theta^2 \wedge \theta^4 \\
  d\theta^2=-\theta^1 \wedge \alpha_2-\theta^2 \wedge \alpha_3-\theta^3 \wedge \theta^4 \\
  d\theta^3=-\theta^2 \wedge \alpha_2+\theta^3 \wedge \alpha_1-2\theta^3 \wedge \alpha_3\\
  d\theta^4=-\theta^4 \wedge \alpha_1+\theta^4 \wedge \alpha_3 \\
  d\alpha_1=T^5_{14}~\theta^1 \wedge \theta^4+\theta^4 \wedge \alpha_2 \\
  d\alpha_2=\frac{1}{2}T^5_{14}~\theta^2 \wedge \theta^4-\alpha_1 \wedge \alpha_2 -\alpha_2 \wedge \alpha_3\\
  d\alpha_3=\frac{1}{2}T^5_{14}~\theta^1 \wedge \theta^4 \\
\end{array}
\end{equation}
The invariant structure function for the prolonged coframe is
\begin{equation}\label{bb35}
\begin{array}{ll}
T^5_{14}=&-\frac{2}{9}\frac{a_3 I_4}{a_1^2},\\
\end{array}
\end{equation}
where
\begin{equation}\label{bbb35}
\begin{array}{l}
  I_4=f_{q,q} \left( {f_q ^2+9\,f_p  - 3\,D_x f_q}\right) - 9\,f_{p,p}  + 18\,f_{u,q}  - 6\,f_q f_{p,q}.\\
  \end{array}
\end{equation}
Similarly, the relative invariant $I_4=0$  for the canonical form
$u''' = 0$. Thus, we choose the following branch.
\section*{Branch 1.1. $I_4=0$.}
In this branch, the invariant structure of the prolonged coframe
are all constant. We have produced an invariant coframe with rank
zero on the seven-dimensional space coordinates $x,u,p, q, a_1,
a_2, a_3$. Any such differential equation admits a
seven-dimensional symmetry group of point transformations.

Moreover, the symmetrical version of the Cartan formulation
$\textrm{mod}~(\omega^1,\omega^2,\omega^3)$ is
\begin{equation}\label{b37}
\begin{array}{l}
\frac{a_1}{a_3}dx=\frac{\bar{a}_1}{\bar{a}_3}d\bar{x}\\
\frac{da_1}{a_1}-\frac{a_2}{a_3}dx=\frac{d\bar{a}_1}{\bar{a}_1}-\frac{\bar{a}_2}{\bar{a}_3}d\bar{x},\\
\frac{da_2}{a_1}-\frac{a_2da_3}{a_1a_3}+\frac{1}{18}\left(\frac{a^2_3s_2-9a^2_2-6a_2a_3s_1}{a_1a_3}\right)dx=\frac{d\bar{a}_2}{\bar{a}_1}-\frac{\bar{a}_2d\bar{a}_3}{\bar{a}_1\bar{a}_3}+\frac{1}{18}\left(\frac{\bar{a}^2_3\bar{s}_2-9\bar{a}^2_2-6\bar{a}_2\bar{a}_3\bar{s}_1}{\bar{a}_1\bar{a}_3}\right)d\bar{x},\\
\frac{da_3}{a_3}+\frac{1}{3}s_1~dx=\frac{d\bar{a}_3}{\bar{a}_3}+\frac{1}{3}\bar{s}_1d\bar{x}.\\
\end{array}
\end{equation}
Inserting the point transformation (\ref{ccc}) into (\ref{b37})
and using  $\bar{s}_1=\bar{s}_2=0$ for $\bar{f}=0$ and
$\bar{a_1}=1,\bar{a_2}=0,\bar{a_3}=1$, results in
\begin{equation}\label{b39}
\begin{array}{l}
\frac{a_1}{a_3}=D_x\phi\\
\frac{D_x a_1}{a_1}=\frac{a_2}{a_3},\\
\frac{D_xa_2}{a_1}-\frac{a_2D_xa_3}{a_1a_3}=-\frac{1}{18}\left(\frac{a^2_3s_2-9a^2_2-6a_2a_3s_1}{a_1a_3}\right),\\
\frac{D_xa_3}{a_3}=-\frac{1}{3}s_1.\\
\end{array}
\end{equation}
This system can be simplified further as
\begin{equation}\label{bb39}
\begin{array}{l}
D_xa_3=-\frac{1}{3} a_3 s_1.\\
D_xa_2=\frac{1}{2}\frac{1}{a_3}a^2_2-\frac{1}{18}a_3s_2,\\
D_x a_1=\frac{a_2}{a_3}a_1,\\
D_x\phi=\frac{a_1}{a_3}\\
\end{array}
\end{equation}
This proves the following theorem.
\begin{theorem}
The necessary and sufficient conditions for equivalence of a scalar
third-order ODE $u''' = f(x,u,u',u'')$ to its canonical form
$\bar{u}'''=0$, with seven point symmetries via {\it point
transformation} (\ref{ccc}), are the identical vanishing of the
relative invariants
\begin{equation}\label{b40}
\begin{array}{l}
  I_1  = f_{q,q,q}  \hfill \\
  I_2  = f_{q,q} ^2  + 6\,f_{p,q,q}  \hfill \\
  I_3  =  {4\,f_q ^3  + 18\,f_q \left( {f_p  - D_x f_q } \right) + 9\,D_x^2 f_q    - 27\,D_x f_p+ 54\,f_u  }\\
  I_4 =f_{q,q} \left( {f_q ^2+9\,f_p  - 3\,D_x f_q}\right) - 9\,f_{p,p}  + 18\,f_{u,q}  - 6\,f_q f_{p,q}.\\
\end{array}
\end{equation}
Given that the the system of relative invariants (\ref{b40}) is
zero, the linearizing point transformation (\ref{ccc}) is defined
by
\begin{equation}\label{b42}
\begin{array}{l}
D_x\phi=\frac{a_1}{a_3},\\
\phi_x\psi_u   - \phi_u \psi_x=\frac{a^2_1}{a_3},\\
\end{array}
\end{equation}
where $a_1(x,u,p), a_2(x,u,p,q),a_3(x,u,p)$ are auxiliary
functions given by
\begin{equation}\label{b43}
\begin{array}{l}
D_xa_3=-\frac{1}{3}f_q a_3,\\
D_xa_2=\frac{1}{2}\frac{1}{a_3}a^2_2-\frac{1}{18}a_3\left( 2f_q ^2 +9\,f_p   - 3\,D_x f_q \right),\\
D_x a_1=\frac{a_2}{a_3}a_1,\\
\left(\frac{a_1}{a_3}\right)_{p,p}=0,\\
\left(\frac{a^2_1}{a_3}\right)_{p}=0.\\
\end{array}
\end{equation}
\end{theorem}
\begin{remark}\label{r1}
The last two equations of the system (\ref{b43}) ensure the
compatibility of the system (\ref{b42}). Moreover, they
provide that
\begin{equation}\label{b44}
\begin{array}{l}
a_2(x,u,p,q)=-\frac{1}{6}a_3 f_{q,q}~q+A(x,u,p),
\end{array}
\end{equation}
for some function $A(x,u,p)$.
\end{remark}
\section{Illustration of the theorem}
\begin{example}\rm \cite{ibr}
Consider the nonlinear ODE
\begin{equation}\label{n1}
\begin{array}{c}
u'''=(\frac{6}{u}u'+\frac{3}{x})u''-\frac{6}{u^2}{u'}^3-\frac{6}{xu}{u'}^2-\frac{6}{x^2}u'-6\frac{u}{x^3}.\\
\end{array}
\end{equation}
The function
\begin{equation}\label{n2}
\begin{array}{ll}
f(x,u,p,q)=(\frac{6}{u}p+\frac{3}{x})q-\frac{6}{u^2}{p}^3-\frac{6}{xu}{p}^2-\frac{6}{x^2}p-6\frac{u}{x^3}\\
\end{array}
\end{equation} satisfies the constraints $I_1=I_2=I_3=I_4=0$;
as a consequence, this equation admits the seven-dimensional point
symmetry group. Moreover, it is equivalent to the canonical form
$\bar{u}'''=0$. We outline the steps.

\begin{description}
\item[Step 1] Solve the system (\ref{b43}) for the auxiliary functions
$a_1(x,u,p), a_2(x,u,p,q),a_3(x,u,p)$ as shown below.

By using the total derivative definition and comparing the
coefficients in the first equation of the system (\ref{b43})
\begin{equation}\label{n3}
\begin{array}{lll}
D_xa_3=-(\frac{2}{u}p+\frac{1}{x})a_3,\\
\end{array}
\end{equation}
yield the system
\begin{equation}\label{n4}
\begin{array}{lll}
\frac{\partial a_3}{\partial p}=0,&\frac{\partial a_3}{\partial u}=-\frac{2}{u}a_3,&\frac{\partial a_3}{\partial x}=-\frac{1}{x}a_3,\\
\end{array}
\end{equation}
which has a solution  $a_3=\frac{1}{x u^2}$.

The second equation of the system (\ref{b43}) reads
\begin{equation}\label{n5}
\begin{array}{lll}
D_xa_2=\frac{1}{2}x u^2 a^2_2+4\,{\frac{{p}^{2}}{x{u}^{4}}}+2\,{\frac {p}{{x}^{2}{u}^{3}}}+\frac{3}{2}\,{ \frac{1}{{x}^{3}{u}^{2}}}-2\,{\frac {q}{x{u}^{3}}},\\
\end{array}
\end{equation}
Using Remark (\ref{r1}), $a_2$ is independent of $q$. Therefore,
comparing the coefficients of the variable $q$ in equation
(\ref{n5}) yield
\begin{equation}\label{n6}
\begin{array}{lll}
\frac{\partial a_2}{\partial p}=-\frac{2}{xu^3}\\
\end{array}
\end{equation}
which gives $a_2=-\frac{2 p}{xu^3}+F_1(x,u)$. Substituting $a_2$
back into equation (\ref{n5}) and then comparing the coefficients of
the variable $p$, result in
\begin{equation}\label{n7}
\begin{array}{lll}
\frac{\partial F_1}{\partial u}=-\frac{2}{u}F_1\\
\end{array}
\end{equation}
which gives $F_1(x,u)=\frac{1}{u^2}F_2(x)$. Substituting
$a_2=-\frac{2p}{xu^3}+\frac{1}{u^2}F_2(x)$ back into equation
(\ref{n5}) yields the Riccati equation
\begin{equation}\label{n8}
\begin{array}{lll}
\frac{d F_2}{d x}=\frac{x}{2}F^2_2+\frac{3}{2x^3}\\
\end{array}
\end{equation}
which gives a solution $F_2(x)=-\frac{1}{x^2}$. Thus,
$a_2=-\frac{2p}{xu^3}-\frac{1}{x^2u^2}$ is a solution of equation
(\ref{n5}).

Now, the third equation of the system (\ref{b43}) reads
\begin{equation}\label{n9}
\begin{array}{l}
D_x a_1=-(\frac{2p}{u}+\frac{1}{x})a_1.\\
\end{array}
\end{equation}
Now by utilizing the total derivative definition and comparing the
coefficients of equation (\ref{n9}), results in
\begin{equation}\label{n10}
\begin{array}{lll}
\frac{\partial a_1}{\partial p}=0,&\frac{\partial a_1}{\partial u}=-\frac{2}{u}a_1,&\frac{\partial a_1}{\partial x}=-\frac{1}{x}a_1,\\
\end{array}
\end{equation}
which has a solution  $a_1=\frac{1}{x u^2}$.

\item[Step 2] Obtain  the linearizing transformation (\ref{ccc}) by
 solving the system (\ref{b42}) as shown below.

The first equation of the system (\ref{b42}) reads
\begin{equation}\label{n11}
\begin{array}{l}
D_x\phi=1,\\
\end{array}
\end{equation}
which has a solution $\phi(x,u)=x$.

Now, the second equation of the system (\ref{b42}) reads
\begin{equation}\label{n12}
\begin{array}{l}
\psi_u =\frac{1}{x u^2},\\
\end{array}
\end{equation}
which has a solution $\psi(x,u)=-\frac{1}{x u}$.

Therefore, the canonical form $\bar{u}'''=0$  can be found for
the ODE (\ref{n1}) via the point transformation
$$\bar{x}=x, \bar{u}=-\frac{1}{xu}.$$

It should be remarked that one can replace $\bar{u}=-\frac{1}{xu}$
with $\bar{u}=\frac{c}{xu}$ for some constant $c$ as
$\bar{u}'''=0$ admits the scaling symmetry.
\end{description}
\end{example}
\begin{example}\rm
We now study the nonlinear ODE
\begin{equation}\label{n13}
\begin{array}{c}
u'''=3\frac{{u''}^2}{1+u'}.\\
\end{array}
\end{equation}
The function
\begin{equation}\label{n14}
\begin{array}{ll}
f(x,u,p,q)=3\frac{{q}^2}{1+p}.\\
\end{array}
\end{equation} satisfies the constraints $I_1=I_2=I_3=I_4=0$;
as a result, this equation admits the seven-dimensional point
symmetry group as well. Moreover, it should be equivalent to the canonical form
$\bar{u}'''=0$. We demonstrate this.

\begin{description}
\item[Step 1] Solve the system (\ref{b43}) for the auxiliary functions
$a_1(x,u,p), a_2(x,u,p,q),a_3(x,u,p)$ as shown below.

By using the total derivative definition and comparing the
coefficients in the first equation of the system (\ref{b43}), one has
\begin{equation}\label{n15}
\begin{array}{lll}
D_xa_3=-\frac{2q}{1+p}a_3,\\
\end{array}
\end{equation}
which has a solution  $a_3=\frac{1}{(1+p)^2}$.

The second equation of the system (\ref{b43}) reads
\begin{equation}\label{n16}
\begin{array}{lll}
D_xa_2=\frac{1}{2}(1+p)^2a^2_2-\frac{1}{2}\frac{q^2}{(1+p)^4}.\\
\end{array}
\end{equation}
Using Remark (\ref{r1}), $a_2=-\frac{q}{(1+p)^3}+A(x,u,p)$, for
some function $A(x,u,p)$. Therefore, substituting $a_2$ into
equation (\ref{n16}) and comparing the coefficients of the
variable $q$, yield
\begin{equation}\label{n17}
\begin{array}{lll}
\frac{\partial A}{\partial p}=-\frac{1}{1+p}A\\
\end{array}
\end{equation}
which gives $A=\frac{F_1(x,u)}{1+p}$. Substituting $a_2$ back into
equation (\ref{n16}) and then comparing the coefficients of the
variable $p$ gives rise to
\begin{equation}\label{n18}
\begin{array}{lll}
\frac{\partial F_1}{\partial x}=\frac{\partial F_1}{\partial u}\\
\end{array}
\end{equation}
which results in $F_1(x,u)=F_2(x+u)$. Substituting
$a_2=-\frac{q}{(1+p)^3}+\frac{F_2(x+u)}{1+p}$ back into equation
(\ref{n16}) yields the Riccati equation
\begin{equation}\label{n19}
\begin{array}{lll}
\dot{F}_2=\frac{1}{2}F^2_2\\
\end{array}
\end{equation}
which has a solution $F_2(x)=0$ and hence
$a_2=-\frac{q}{(1+p)^3}$.

Now, the third equation of the system (\ref{b43}) reads
\begin{equation}\label{n20}
\begin{array}{l}
D_x a_1=-\frac{q}{1+p}a_1.\\
\end{array}
\end{equation}
This has a solution  $a_1=\frac{1}{1+p}$.

\item[Step 2] Obtain  the linearizing transformation (\ref{ccc}) by
 solving the system (\ref{b42}) as shown below.

The first equation of the system (\ref{b42}) reads
\begin{equation}\label{n21}
\begin{array}{l}
D_x\phi=1+p,\\
\end{array}
\end{equation}
which has a solution $\phi(x,u)=x+u$.

Now, the second equation of the system (\ref{b42}) reads
\begin{equation}\label{n22}
\begin{array}{l}
\psi_u-\psi_x =1,\\
\end{array}
\end{equation}
which has a solution $\psi(x,u)=-x$.

Therefore, the canonical form $\bar{u}'''=0$  can be determined
for the ODE (\ref{n13}) via the point transformation
$$\bar{x}=x+u, \bar{u}=-x.$$
\end{description}
\end{example}
\begin{example}\rm
We now focus on the nonlinear ODE
\begin{equation}\label{n9}
\begin{array}{c}
u'''=\frac{3}{2}\frac{{u''}^2}{u'}.\\
\end{array}
\end{equation}
This equation does not satisfy our condition. Indeed the function
\begin{equation}\label{n10}
\begin{array}{ll}
f(x,u,p,q)=\frac{3}{2}\frac{{q}^2}{p}\\
\end{array}
\end{equation} gives $I_2\ne0$. Consequently, this equation cannot admit a seven-dimensional point
symmetry group. This agrees with the fact that the ODE (\ref{n9})
has maximal six point symmetries \cite{IbrFazal1996}.
\end{example}
\section{Conclusion}
There are essentially two basic approaches to linearization via point transformations.
One is algebraic and the second geometric.

In this work, we have successfully applied the Cartan equivalence method
to solve the linearization problem for scalar third-order
ODEs with maximal symmetry Lie algebra to its simplest canonical form, viz.
reduction via invertible map to the simplest linear third-order ODE admitting seven
point symmetries.

What is more, we were able to provide a step by step approach to construct a point
transformation which did the reduction to the simplest canonical form. This has not been
given before using this approach.  We also mention in this context that in
\cite{ibr}, they utilize the direct method, in order to find reduction to the Laguerre-Forsyth
canonical form which can admit four, five or seven point symmetries. Here we stress that we have focused
on the maximal symmetry case, via the Cartan method and have dealt with the construction of the linearizing
maps.

Moreover, we have illustrated our main result by means of examples to demonstrate the effectiveness of our method.

\subsection*{Acknowledgments}
Ahmad Y. Al-Dweik is grateful to the King Fahd University of
Petroleum and Minerals for its continuous support and excellent research
facilities. FMM is indebted to the NRF of South Africa for support.

\end{document}